\documentclass[11pt, twoside]{amsart}
\usepackage{latexsym}
\usepackage{graphicx}
\usepackage[centertags]{amsmath}
\usepackage{amsfonts}
\usepackage{amssymb}
\usepackage{amsthm}
\usepackage{newlfont}
\usepackage{subfigure}
\usepackage{epsfig}

\theoremstyle{plain}
\newtheorem{theorem}{Theorem}[section]

\newtheorem{lemma}[theorem]{Lemma}

\theoremstyle{definition}
\newtheorem{defn}[theorem]{Definition}
\newtheorem{remark}[theorem]{Remark}

\newcommand{\s}[1][1]{\mathbf{S}^{#1}}

\newcommand{\cross}{\times}
\newcommand{\bndry}{\partial}

\newcommand{\pii}{\pi_1}

\newcommand{\tc}{\cal T}

\begin{document}
\title[The word problem and injective handlebodies]{The word problem for 3-manifolds built from injective handlebodies.}

\author{James Coffey}
\address{University of Melbourne}
\email{coffey@ms.unimelb.edu.au}
\thanks{Research partially supported by the Australian Research Council}

\abstract{
    This paper gives a proof that the fundamental group of a  class of closed orientable 3-manifolds constructed from three injective handlebodies has a solvable word problem.  This is done by giving an algorithm to decide if a closed curve in the manifold is null-homotopic. Non-Haken and non-Seifert fibered examples are constructed by performing Dehn surgery on a class of two-bridge knots.}}

\maketitle


\section{Introduction}
The word problem, which was first stated by M. Dehn in \cite{de1}, along with the conjugacy problem and the isomorphism problem lie at the core of combinatorial group theory. The word problem asks that; given a group with finite presentation, is it possible, given any  word in its generators, to determine in a finite time if it  represents the trivial element. In general, a group is shown to have  a solvable word problem by giving an algorithm that determines if a given word is trivial and terminates in finite time.  There are many examples of groups that do not have a solvable word problem. However, it is a long standing conjecture that the fundamental group of a 3-manifold will always have solvable word problem.

Even though the word problem is an algebraic question, when we are looking at it in relation to the  fundamental group of manifolds, the question can by restated geometrically.  As any word in the fundamental group of a manifold can be represented be a loop in the manifold the word problem becomes; given an immersion of $\s$ in some manifold it is possible to determine if the map is null-homotopic.  Thus to show that fundamental group of a 3-manifold has solvable word problem an algorithm is given, which terminates in finite time and determines if a given loop bounds an immersed disk.

F. Waldhausen in \cite{wa2} showed that the fundamental group of all $P^2$-irreducible Haken 3-manifolds has a solvable word problem. As having a solvable word problem is a property inherited by a sub-groups, all virtually Haken manifolds also have solvable word problem. Andrew Skinner showed in \cite{sk1} that the fundamental group of 3-manifolds has a solvable word problem if the 3-manifold contains an essential singular surface that is not $\s[2]$ or $P^2$ and satisfies the 1-line and 1-point property.  This condition was defined by J.Hass and P.Scott in \cite{hs1} and essentially means that the surfaces self intersections are `nice'. If Thurston's Geometrisation conjecture is proven, which looks likely following the announcements of Perelman, then the conjecture will be proven as a result.  However, it is still an interesting question to answer, using more direct means, for topological classes of manifolds that do not appear to have any obvious unifying underlying geometry. Also understanding why 3-manifold groups have solvable word problem can give hints to other classes of groups with this property.

The class of manifolds being considered in this paper are said to meet the `disk-condition'.  The construction of these manifolds was introduced by H. Rubinstein and the author in \cite{C&R1}.  This construction is an extension to handlebodies of the gluing of three solid tori which produces non-Haken Seifert fibered manifolds with infinite fundamental group. A definition is given in the next section. The main theorem proved in this paper is:

\begin{theorem}
    If $M$ is a 3-manifold that satisfies the disk-condition then the word problem in $\pii(M)$ is solvable.
\end{theorem}

This theorem is proven by giving an algorithm to test if a given word in $\pii(M)$ is trivial and showing that the algorithm  terminates in bounded time.   The idea is that if a loop is trivial that it must bound some singular disk. Therefore the algorithm tries to find a disk bounded by the loop by progressively filling it in one face at a time.  There are two steps that contain the majority of the work to show there  is an algorithm.  The first is to show that if a loop is trivial there is always a face with certain properties that makes up part of  the disk and secondly, that the face can indeed be found.  This proof is quite satisfying as, once the structure of the manifold is set up, the result is obtained directly and easily.

To help motivate this paper, in section \ref{section: non-haken example}, it is shown that Dehn surgery on a family of two-bridge knots produces non-Haken  non-Seifert fibered manifolds which meet the disk-condition.

\section{Preliminaries and definitions}

Throughout this paper, we will assume that, unless stated otherwise, we are working in the PL category of manifolds and maps.  Even though we will not explicitly use this structure we will use ideas that are a consequence, such as regular neighbourhoods and transversality as defined by C. Rourke and B. Sanderson in \cite{rou&sa1}. The standard definitions in this field, as given by J. Hempel in \cite{Hem1} or W. Jaco in \cite{Ja1}, are used.

A manifold $M$ is \textbf{closed} if it is compact and $\bndry M = \emptyset$ and \textbf{irreducible} if every embedded $S^2$ bounds a ball.  We will assume, unless otherwise stated that all 3-manifolds are orientable. The reason for this is that all closed non-orientable $\mathbb P^2$-irreducible 3-manifolds are Haken. (A manifold is $\mathbb P^2$-irreducible  if it is irreducible and does not contain any embedded $2$-sided projective planes).

If $M$ is a 3-manifold and $S$ is some surface, which is not a sphere, disk or projective plane,  the map $f:S\to M$ is called \textbf{essential} if the induced map $f_*:\pii(S)\to \pii(M)$ is injective. This is also known as a $\mathbf{\pi_1}$\textbf{-injective} map.   Also $f:S \to M$ is a \textbf{proper map} if $f(S)\cap \bndry M  = f(\bndry S)$.  If $F:S\cross I \to M$ is a homotopy/isotopy such that $F(S, 0)$ is a proper map, then it assumed, unless otherwise stated  that $F(S,t)$ is a proper map for all $t\in I$.  To reduce  notation, an isotopy/homotopy of a surface $S\subset M$ is used without defining the map.  Here we are assuming that there is a map $f:S \to M$ and we are referring to  an isotopy/homotopy of $f$. Defining the map is often unnecessary and would only add to excessive book keeping.

If $H$ is a handlebody and $D$ is a properly embedded disk in $H$ such that $\bndry D$ is essential in $\bndry H$ then $D$ is a  \textbf{meridian disk} of $H$.  If $D$ is a proper singular disk in $H$ such that $\bndry D$ is essential in $\bndry H$, then it is called a \textbf{singular meridian disk}.

In this paper normal curve theory, as defined by S. Matveev in \cite{Mat1}, is used to list finite classes of curves in surfaces.  This definition uses a triangulation of the surface to define normal curves. The surfaces in this paper have polygonal faces, but a barycentric subdivision will produce the required triangulation.   CAT$(0)$ metric spaces, as defined by M.Bridson and A.Haefliger in \cite{b&h}, are also used in this paper.

\begin{defn}
    For $H$ a handlebody, $\tc$ a set of curves in $\bndry H$ and $D$ a meridian disk, let $|D|$ be the number of intersection between $D$ and $\tc$.
\end{defn}

Freedman, Hass and Scott showed in \cite{fhs1} that if we put a metric on $\bndry H$ and isotop both $\tc$ and $D$ so that $\tc$ and $\bndry D$ are length minimising in $\bndry H$ then the number of intersections will be minimal. Note that when there are parallel curves in $\tc$ we need to `flatten' the metric in their neighbourhood so they remain disjoint.  Let $\mathcal D$ be a maximal set, up to ambient isotopy, of meridian disks for $H$, where each element of $\mathcal D$ is represented by a meridian disk whose boundary is a length minimising geodesic.  Therefore the number of intersections between the boundaries of two disks in $\cal D$ is minimal and in  \cite{C&R1} the following lemma is proven.

\begin{lemma}\label{lemma: m-disk intersection.}
    Any two disks of $\mathcal D$ can be isotoped, leaving their boundaries fixed, so that the curves of intersections are  properly embedded arcs.
\end{lemma}

\begin{remark}
This lemma is proven by constructing an isotopy to remove curves of intersection which are closed.  This was done using the usual innermost argument and the fact that a handlebody is irreducible.
\end{remark}

We will assume from this point on that all curves in $\bndry H$ and meridian disks have been isotoped to have minimal intersection.

\begin{defn}
    If $H$ is a handlebody and $\tc$ is a set of essential disjoint simple closed curves in $\bndry H$, then $\tc$ meets the $\mathbf n$ \textbf{disk-condition} in $H$ if for every meridian disk $D$, $|D|\geq n$.
\end{defn}

See \cite{C&R1} for sufficient and necessary condition for $\tc$ to satisfy the $n$ disk-condition in $H$.  It is also shown in \cite{C&R1} that there is an algorithm to decide if a set of curves in $\bndry H$ meets the $n$ disk-condition.

Next we are going to give a description of the construction of 3-manifolds that meet the `disk-condition'.  Let $H_1$, $H_2$ and $H_3$ be three handlebodies.  Let $S_{i,j}$, for $i\not= j$, be a sub-surface of $\bndry H_i$ such that:

\begin{enumerate}
    \item
        $\bndry S_{i,j} \not= \emptyset$,
    \item
        the induced map of $\pii(S_{i,j})$ into $\pii(H_i)$ is injective,
    \item
        for $j\not= k$, $S_{i,j}\cup S_{i,k} = \bndry H_i$,
    \item
        $\tc_i  = S_{i,j}\cap S_{i,k} = \bndry S_{i,j} = \bndry S_{i,k}$ is a set of disjoint essential simple closed curves that meet the $n_i$ disk-condition in $H_i$,
    \item
        and  $S_{i,j}\subset \bndry H_i$ is homeomorphic to $S_{j,i} \subset \bndry H_j$.
\end{enumerate}

\begin{figure}[h]
   \begin{center}
       \includegraphics[width=6cm]{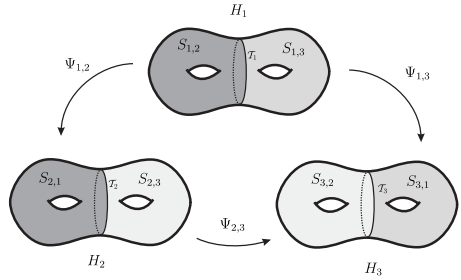}
    \end{center}
   \caption{Homeomorphisms between boundaries of handlebodies.} \label{fig:handlebodies}
\end{figure}

Note that $S_i$ does not need to be connected and that the handlebodies can be different genus. Now that we have the boundary of each handlebody cut up into essential faces  we want to glue them together by homeomorphisms, $\Psi_{i,j}: S_{i,j}\to S_{j,i}$, that agree along $\tc_i$'s, as in figure \ref{fig:handlebodies}.  The result is a closed 3-manifold $M$, for which the image of each handlebody is embedded.

\begin{defn}
    If $M$ is a manifold constructed from three handlebodies as above such that $\tc_i$ meets the $n_i$ disk-condition in $H_i$ and

    \begin{equation}\label{equation: disk-condition}
         \sum_{i=1,2,3} \frac{1}{n_i} \leq \frac{1}{2}
    \end{equation}

    then $M$ meets the $\mathbf{(n_1,n_2,n_3)}$ \textbf{disk-condition}. If we are not talking about a specific $(n_1,n_2,n_3)$, the manifold is said to meet just the \textbf{disk-condition}.
\end{defn}

Although the above definition of the construction of manifolds is technically correct for the purposes of this paper, it is not the best way to view them. Thus from this point on, we will view 3-manifolds that meet the disk-condition in the following way.  Assume that $M$ is a manifold that meets the disk-condition and $H_1$, $H_2$ and $H_3$ are the images of the handlebodies of the previous definition in $M$.  Then $M = \bigcup_{i=1,2,3} H_i$ and each $H_i$ is embedded in $M$. Then $X = \bigcup_{i=1,2,3} \bndry H_i$ is a 2-complex that cuts $M$ up into handlebodies. Note that $X$ does not meet the usual definition of 2-complex, as it is constructed by gluing surfaces along their boundaries. These surfaces can easily be cut up into cells.  Also $\tc = \bigcap_{i=1,2,3} H_i$ is a set of essential disjoint simple closed curves in $M$ that meet the $n_i$ disk-condition in $H_i$ where $\sum_{i=1,2,3} 1/n_i\leq 1/2$.

The following lemma about the homotoping of maps of disks into manifolds meeting the disk-condition is proven in \cite{C&R1}.

\begin{lemma}\label{lemma: removing pullback graphs}
    Let $M$ be a manifold that meets the disk condition and $D$ be a disk. If $f:D\to M$ is a map such that $f(\bndry D)\subset int(H_i)$, for some $i$, then $f$ can be homotoped so that $f(D)\subset int(H_i)$.
\end{lemma}

This lemma is proven by using the usual Euler characteristic arguments and the disk-condition to construct the homotopy of $f$. If the image of the disk intersects $X= \bigcup H_i$'s transversely, then $\Gamma = f^{-1}(X)$ is a trivalent graph that cuts $D$ up into faces.  By the Euler characteristic and the disk-condition, at least one of the faces, which is a disk, maps to a disk parallel into a boundary of a handlebody. This gives a homotopy of $f$ to remove the face and reduce $\Gamma$.

\begin{defn}
    Let $H$ be a genus $g$ handlebody. We shall call $\mathbf D \subset \mathcal D$ a \textbf{system of meridian disks} if all the disks are disjoint, non parallel and they cut $H$ up into a set of 3-balls. If $\mathbf D$ cuts $\bndry H$ up into $2g-2$ pairs of pants (thrice punctured 2-spheres) then it is a \textbf{basis} for $H$.
\end{defn}

If $H$ has genus $g$, then a minimal system of meridian disks for $H$ consists of $g$ disks and it cuts $H$ up into a single ball.

\begin{defn}
    Let $P$ be a punctured sphere and $\gamma$ be a properly embedded arc in $P$. If  both ends of $\gamma$ are in the one boundary component of $\bndry P$ then it is called a \textbf{wave}. If a wave in $P$ is not isotopic into $\bndry P$ then it may be called an \textbf{essential wave}.
\end{defn}

Let $H$ be a handlebody, $\tc$ a set of essential disjoint simple closed curves in $\bndry H$, $\mathbf D$ be a system (basis) of meridian disks for $H$ and $\tc\subset \bndry H$ a set of essential simple closed curves. Let  $\{ P_1,...,P_l\}$ be the resulting set of punctured spheres produced when we cut $\bndry H$ along $\mathbf D$ and $\tc_i = P_i \cap \tc$.  Thus $\tc_i$ is a set of properly embedded disjoint arcs in $P_i$.  Note that as we are assuming the intersections between $\tc$ and $\mathbf D$ are minimal, all waves in the $\tc_i$'s are essential.

\begin{defn}
    If each $\tc_i$ contains no waves then $\mathbf D$ is said to be a \textbf{waveless} system (basis) of meridian disks for $H$.
\end{defn}

\begin{defn}\label{defn: n-waveless}
    Let $\mathbf D$ be a waveless system of disks, if for $1\leq i \leq l$ any essential wave $\gamma \subset P_i$, $|\gamma \cap \tc_i|\geq n /2$, then $\mathbf D$ is called an $\mathbf n$\textbf{-waveless} system (basis) of meridian disks.
\end{defn}

The following lemma was proven in \cite{C&R1}.

\begin{lemma}\label{lemma: finding n-waveless system}
    If $H$ is a handlebody and $\tc \subset \bndry H$ is a set of essential curves that meet the $n$ disk-condition, then there an algorithm to find an $n$-waveless minimal system of meridian disks.
\end{lemma}

\section{Non-Haken example} \label{section: non-haken example}

The following lemma gives a sufficient condition for a set of curves in the boundary  of handlebody to meet the $n$ disk-condition.

\begin{lemma}\label{lemma: sufficient condition}
Let $H$ be a handlebody and $\tc$ be a set of pairwise disjoint essential simple closed curves in $\bndry H$. If $H$ has an $n$-waveless basis $\mathbf D$ then $\tc$ meets the $n$ disk-condition in $H$.
\end{lemma}

\begin{proof}
The full proof of this lemma is in \cite{C&R1}.  By the definition of $n$-waveless basis, any meridian  disk $D$ properly isotopic to disk in $\mathbf D$ must intersect $\tc$ at least $n$ times.  If $D$ is a meridian disk not isotopic to a disk in $\mathbf D$, then $D$ must have an intersection with at least one disk in $\mathbf D$ that can not be removed by isotopy. Therefore $\bndry D$ contains at least two essential waves in $\overline{\bndry H - \mathbf D}$ and thus, by the definition of $n$-waveless basis, $D$ intersects $\tc$ at least $n$ times.
\end{proof}

\begin{remark}
In \cite{C&R1} it is shown that $H$ having a minimal $n$-waveless system of meridian disks for $\tc$ is a necessary and sufficient condition for $\tc$ to meet the $n$ disk-condition in $H$.
\end{remark}

\begin{figure}[h]
   \begin{center}
       \includegraphics[width=4cm]{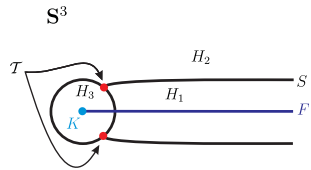}
    \end{center}
   \caption{Handlebodies in Dehn filling construction.}
   \label{fig: Handlebodies in Dehn filling construction.}
\end{figure}

Let $K$ be a knot in $\s[3]$ and $F$ be a free spanning surface such that $K$ meets the $3$ disk-condition in the handlebody $\overline{\s[3] - F}$. If we perform $(p,q)$ Dehn surgery along $K$ where $|p|\geq 6$, then the resulting manifold $M$ will meet the disk-condition. To prove this we need to construct a 2-complex $X \subset M$ which cuts $M$ up into injective handlebodies.  Let $n(K)$ be a regular neighbourhood of $K$ and $n(F)$ a regular neighbourhood of $F$. Let $H_1 =  \overline{ n(F)- n(K)}$, $H_2 = \overline{\s[3] - (n(K)\cup H_1)}$, $X = \bndry H_1 \cup \bndry H_2$ and $\tc = H_1\cap H_2 \cap n(K)$, as in figure \ref{fig: Handlebodies in Dehn filling construction.}. Therefore $\tc$ is two copies of $K$.  We can fiber $H_1$ as an $I$-bundle over $F' = \overline{F - n(K)}$ where $\tc$ is the boundary curves of the vertical boundary. The horizontal surface is either two copies of $F'$ or it double covers $F'$, depending on whether $F'$ is orientable or not. If $\gamma$ is properly embedded arc in $F'$ which is not isotopic into $\bndry F'$, then the lift of $\gamma$ to $H_1$ is a meridian disk. If $\Gamma$ is a set of properly embedded arcs whose lift to $H_1$ is a basis $\mathbf D$, then $\overline{F'-\Gamma}$ is a set of disks. By the $I$-bundle structure, any essential wave in $\overline{\bndry H_1-\mathbf D}$ must intersect $\tc$ at least twice. Therefore $\mathbf D$ is a 4-waveless basis for $H_1$ and  by lemma \ref{lemma: sufficient condition} $\tc$ meets the 4 disk-condition in $H_1$. The manifold $\overline{\s[3] - (F \cup H_2)}$ is homeomorphic to $\bndry H_2 \cross I$. As $K$ meets the  $6$ disk-condition in $\overline{\s[3] - F}$ and $\tc$ is  two copies of $K$, $\tc$ meets the $6$ disk-condition in $H_2$.  Remove $int(n(K))$ to produce a manifold $M'$. Then glue a solid torus $H_3$ back onto $\bndry M'$ such that the meridian disk of $H_3$ winds meridainly at least 6 times around $\bndry M'$ and  thus intersects $\tc$ at least 12 times. The resulting manifold will meet the $(4,6,12)$ disk-condition.

\begin{figure}[h]
   \begin{center}
       \includegraphics[width=7cm]{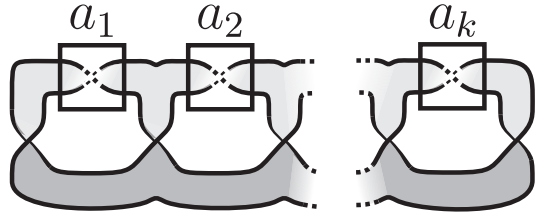}
    \end{center}
   \caption{`Plumbed' twisted bands.}
   \label{fig: plumbed twisted bands.}
\end{figure}

Any two-bridge (rational) link $K$ has a projection that is the boundary of $k$ twisted  twisted bands `plumbed' together, as shown in figure \ref{fig: plumbed twisted bands.} where the $i$'th band has $a_i$ right handed twists.  The arcs in $F$ (figure \ref{fig: plumbed twisted bands.}) can be used like the component arcs in normal curves to construct the boundary curves of the meridian disks for a basis of $\overline{\s[3] - F}$. If $K$  is a knot and $|a_i|\geq 3$ for $1\geq i \geq k$,  then the basis constructed is a 3-waveless basis. Therefore from  lemma \ref{lemma: sufficient condition} we know  that  $K$ meets the 6 disk-condition in $\overline{\s[3] - F}$.  From above if we take $(p,q)$ Dehn surgery, such that $|p| \geq 6$, then the resulting manifold meets the $(4,6,12)$ disk-condition. Hatcher and Thurston showed in \cite{H&T1} that all but finitely many of the manifolds resulting from Dehn surgeries on two-bridge knots are non-Haken and can not be Seifert fibered. Further examples of manifolds meeting the disk-condition are constructed using Dehn surgery and branched covers in \cite{C&R1}.

\section{The word problem}

In this section we give a proof that the word problem is solvable in the fundamental group of manifolds meeting the disk-condition.  In the first subsection some necessary definitions are given and then the algorithm for solving the word problem is outlined.  In the following subsection the lemmas need to prove that the algorithm indeed works are given and proved.

\subsection{The algorithm.}
\label{section: the algoirthm}

Let $M$ be a 3-manifold that meets the disk-condition. As described earlier, let $H_i\subset M$ for $1\geq i \geq 3$, be the three injective handlebodies such that $M= \bigcup H_i$. Then $X=\bigcup \bndry H_i$ is a 2-complex and $\tc = \bigcap H_i$ is the triple curves.  Let $D$ be  a disk and the $f:D\to M$ be transverse to $X$.  Then $\Gamma = f^{-1}(X)$ is a trivalent graph that cuts $D$ up into faces.  In \cite{C&R1} it is shown that there is a homotopy of $f$, leaving the boundary fixed, such that all the faces of $D$ are disks.  This is proven using the usual Euler characteristic arguments.  A face $A$ of $D$ is called a boundary face if $A \cap \bndry D \ne \emptyset$.

Let $\alpha$ be a simple closed loop in $M$. Note that saying $\alpha$ is null-homotopic and that it bounds a disk are equivalent statements. The general idea of the algorithm is to try and construct a homotopy of $\alpha$ to a point by looking along $\alpha$ for possible boundary faces.  Then $\alpha$ is homotoped  across each of these faces in turn to produce a number of new loops.  The process is then repeated for all the new curves. This process can also be thought of as constructing a disk bounded by $\alpha$. First we need to define the boundary faces that we are looking for:

\begin{defn}
    Let  $H$ be a handlebody and $\tc$ a set of disjoint curves in its boundary that meet the $n$ disk-condition.  Let $\beta$ be a properly embedded arc in $H$ disjoint to $\tc$. Let $\beta'\subset \bndry H$ be an arc homotopic to $\beta$, keeping its boundary fixed.  If $\beta'$ intersects $\tc$ less than $n/2$ times then it is a \textbf{short boundary path} of $\beta$. $\beta$ is \textbf{short} if it has an associated short boundary path.  If $D$ is a disk such that $\bndry D = \beta \cup \beta'$ and $D \cap \bndry H = \beta'$, then it is called a \textbf{short disk} associated with $\beta$.
\end{defn}

\begin{figure}[h]
  \begin{center}
       \includegraphics[width=6cm]{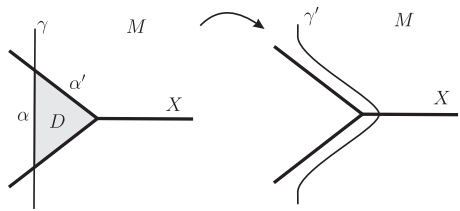}
   \end{center}
  \caption{Short disk move.}
  \label{fig: removing a short disk.}
\end{figure}

Let $f:\s[1] \to M$ be  a map such that $\alpha = f(\s)$ is a simple closed curve transverse to $X$. If we assume that $X\cap \alpha \ne \emptyset$, then $X$ cuts $\alpha$ up into sub-arcs $\alpha_1,...,\alpha_l$. For $\alpha_j$ a short arc in $H_i$ there is an associated short boundary path $\alpha_j' \subset \bndry H_i$ and short disk $D$. $D$ defines a homotopy of $f$ that produces a new curve $\alpha'$.  The homotopy is the identity outside a neighbourhood of $\alpha_j$ and pushes $\alpha$ across $D$ to $\alpha_j' $ and then out the other side of $\bndry H_i$, as shown in figure \ref{fig: removing a short disk.}. This will be referred to as a \textbf{short disk  move}.

Let $f:\s\to M$ such that $\alpha = f(\s)$ is a simple closed curve transverse to $X$.  Let $X$ cut $\alpha$ into sub arcs ${\alpha_1,... \alpha_l}$.  Then the algorithm to solve the word problem in a manifold $M$ that meets the disk-condition is:

\begin{enumerate}
    \item check ${\alpha_j}$'s for short curves and for each short curve find an associated short boundary path,

    \item perform all the possible short disk moves to produce a new set of curves (one for each possible short disk move),

    \item check all the new curves to see if any are disjoint from $X$.  If one is found then terminate, otherwise repeat the process on all the new curves.
\end{enumerate}

A loop disjoint to $X$ will be in some $H_i$ and as they are $\pii$-injective the loop is easily checked to be null-homotopic.

For this algorithm to work we need to show that:
\begin{itemize}
    \item if $\alpha$ is null-homotopic then $\{\alpha_j\}$ always contains at least one short path,

    \item each short path $\alpha_j \subset H_i$ has a finite number of associated short boundary paths up to homotopy in $\bndry H_i$ keeping its ends fixed,

    \item we can always recognise  the short paths in $\{\alpha_i\}$ and find their associated short boundary paths,

    \item and finally, that we know when to stop looking. That is, we have a bound on the number of faces of a disk bound by our original curve $\alpha$ in terms of $l$,  where $l$ is a number of intersections that $l$ has with $X$.
\end{itemize}

\subsection{Short boundary curves.}

If $\alpha$ is null-homotopic then short $\alpha_i$'s and their short boundary path define a face of a disk bounded by a loop in $M$. Therefore we need to show first that if a loop is homotopically trivial that it must have some short sub-arcs.

\begin{lemma}
    If $\alpha$ is a null-homotopic loop in $M$, then at least four of the $\alpha_j$'s must be short.
\end{lemma}

\begin{proof}
Let $D$ be a disk. If $\alpha$ is a homotopically trivial map then there is a map $f:D\to M$ such that $f(\bndry D) = \alpha$. Let $f(D) = D'$ and $\Gamma = \bndry D \cup f^{-1}(X)$.  Therefore $\Gamma$ is a trivalent graph that cuts $D$ up in to faces. As usual we will let $\Gamma$ cut $D$ up into faces. By lemma \ref{lemma: removing pullback graphs} there is a homotopy of $f$  to remove all components of $\Gamma$ contained in $int(D)$.  Thus $\Gamma$ is a connected graph such that $\bndry D \subset \Gamma$ and each face $F_j$ of $D$ is an $(m,n)$-gon. An $(m,n)$-gon is a disk with $m$ vertices in its boundary and that gets mapped by $f$ into a handlebody that meets the $n$ disk-condition.

Next we want to put a metric on $D$.   Let $\cal I$ be the set of internal faces of $D$, that is faces that are disjoint from $\bndry D$ and let $\cal B$ be the set of boundary faces, that is faces that have an edge in $\bndry D$.  We will let each edge be a geodesic arc.  For a face $F\in \cal I$ that is an $(m,n)$-gon let the internal angle at each vertex in its boundary be $\pi - 2\pi/n$ and allow a cone point in $int(F)$, as shown in figure \ref{fig: Angles on faces.}.a. If $K(F)$ is the curvature at the cone point, then by the Gauss-Bonnet theorem we know that  $K(F) = 2\pi(1- m/n)$.  This means that if $m<n$ then $K(F)>0$ and if $m\geq n$ then $K(F)\leq 0$. For the time being we will assume that there are no $(2,n)$-gons in $\cal B$.  At the end of the proof we will show this does not affect the result. For $F \in \cal B$ let the angle at the vertices in the boundary of $F$ disjoint from $\bndry D$ be $\pi - 2\pi/n$.  Note that this means if  $M$ meets one of the minimal disk-conditions, ie $(6,6,6)$, then the sum of the angles around each internal vertex is $2\pi$.  For the vertices $F$ in $\bndry D$ let the internal angle be $\pi/2$, as shown in figure \ref{fig: Angles on faces.}.b.  This agrees with the picture of $f(D)$ in $M$ as $\alpha$ is transverse to $X$.   This means that the total angle at each vertex in $\bndry D$ is $\pi$ and thus $\bndry D$ is geodesic.  All faces of $\cal B$ have at least two vertices in $\bndry D$.   By the Gauss-Bonnet theorem if $F\in \cal B$ is an $(m,n)$-gon and has exactly two vertices in $\bndry D$, then $K(F) = 2\pi - (\pi +2\pi/n(m-2) = \pi - 2\pi/n(m-2)$.  Thus if $F$ has two vertices in $\bndry D$ and $m<n/2+1$ then $K(F) >0$, otherwise  $K(F) \leq 0$. If $F \subset \cal B$ has more than two vertices in $\bndry D$ then it must have at least four and $K(F)\leq 0 $.

\begin{figure}[h]
\begin{center}
  \subfigure[Internal $(m,n)$-gon]{
    \includegraphics[width=4cm]{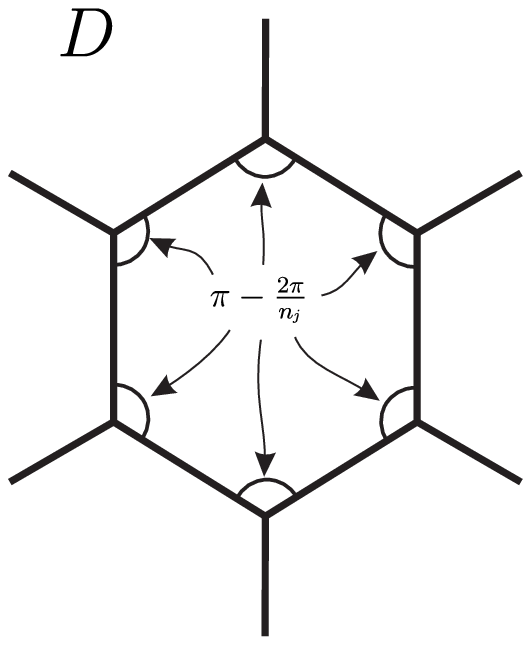}}
    \label{fig: Angles on faces a.}
  \hspace{2cm}
  \subfigure[Boundary $(m,n)$-gon.]{
    \includegraphics[width=4cm]{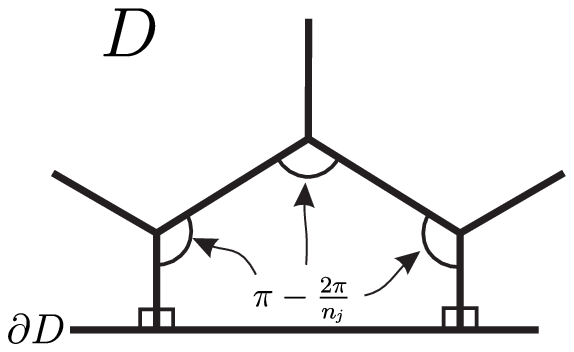}}
    \label{fig: Angles on faces b.}
\end{center}
\caption{Angles on faces of $\Gamma$.}
\label{fig: Angles on faces.}
\end{figure}

We can homotop $f$ so that all faces in $\cal I$ are $(m,n)$-gons such that $m\geq n$. If any faces of $\cal I$ are $(m,n)$-gons where $m<n$, then by the disk-condition the image of such a face, in some $H_i$, is parallel into the boundary.  Therefore there is a homotopy of $f$ to remove the face from $\Gamma$.  As a result the number of faces in $\cal I$ goes down, so this process must come to an end.  We can then assume that each face is an $(m,n)$-gon such that $m\geq n$ and thus $\sum_{\cal I} K(F) \leq 0$.

By assumption all the faces in $\cal B$ are $(m,n)$-gons for $m \geq 3$ and $n\geq 4$ and therefore for any $F\in \cal B$, $K(F) \leq\pi/2$. As $\bndry D$ is geodesic and by Gauss-Bonnet theorem $K(D) = 2\pi \leq \sum_{\cal I}K(F) + \sum_{\cal B}K(F) \leq \sum_{\cal B}K(F)$.  Thus $\cal B$ must contain at least four faces such that they have exactly two vertices in $\bndry D$ and $m \leq n/2+1$.  Let $F \in \cal B$ be an $(m,n)$-gon mapped in $H_i$ such that $m \leq n/2+1$ and it has two vertices in $\bndry D$. Then  $\alpha_j = \overline{f(\bndry F) \cap int(H_i)}$ is a short arc as  $\alpha_j' = f(F) \cap \bndry H_i $ is an associated short boundary path.

Finally we need to show that $(2,n)$-gons can be disregarded in this proof. Let $F\in \cal B$ be an $(m,n)$-gon such that $m \geq 3$. Then if we add a $(2,m)$-gon $F'$, as shown in figure \ref{fig: adding 2-gon to boundary face.}, the total curvature of the two faces is unchanged.  That is, $K(F) = K(F')+K(F'')$.  Therefore, even though the curvature of a $(2,n)$-gon is $\pi$, $\cal B$ must contain at least four $(m,n)$-gons such that $2< m \leq n/2+1$  and have two vertices in $\bndry D$.
\end{proof}

\begin{figure}[h]
  \begin{center}
      \includegraphics[width=8cm]{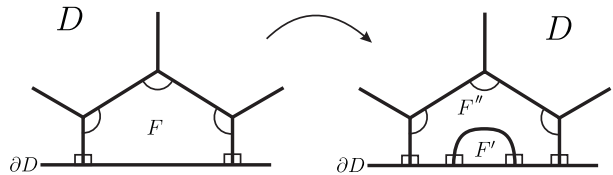}
   \end{center}
  \caption{Adding a $(2,n_j)$-gon to boundary face.}
  \label{fig: adding 2-gon to a boundary face.}
\end{figure}

\begin{lemma}
    Let  $H$ be a handlebody and $\tc$ a set of disjoint curves in its boundary that meet the $n$ disk-condition.  Let $\beta$ be a properly embedded arc in $H$ disjoint to $\tc$. Then $\beta$ has at most one short boundary path up to homotopy in $\bndry H$ leaving the endpoints fixed.
\end{lemma}

\begin{proof}
Assume that the above lemma is incorrect and $\beta$ has two associated short boundary paths $\beta'$ and $\beta''$ that are not homotopic in $\bndry H$ relative to their boundary.  We know that $\beta \cup \beta'$ bounds a short disk $D'$ and $\beta \cup \beta''$ bounds a short disk $D''$. Thus $D' \cap D'' = \beta$ and $D = D'\cup D''$ is a singular disk such that $D\cap \bndry H = \bndry D$.  However $D$ must intersect $\tc$ less than $n$ times. Therefore by the disk-condition there is a proper homotopy of $D$ into $\bndry H$.
\end{proof}

Next we need to show that given a properly embedded arc in a handlebody that we can indeed work out if it is short and, if it is, find an associated short boundary path.

\begin{lemma}\label{lemma: finding short arcs}
    Let $H$ be a handlebody, $\tc$ be a set of simple closed essential arcs in $\bndry H$ that meet the $n$ disk-condition and $\alpha$ be a properly embedded curve in $H$ disjoint to $\tc$. Then there is an algorithm to determine whether $\alpha$ is short and if it is to find an associated short boundary path.
\end{lemma}

\begin{figure}[h]
  \begin{center}
      \includegraphics[width=7cm]{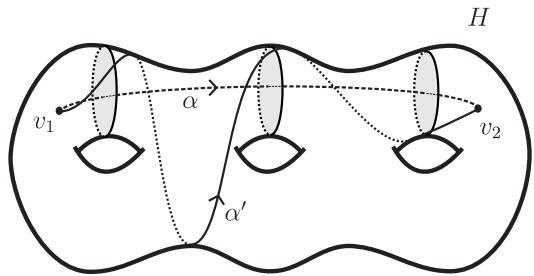}
   \end{center}
  \caption{A short path and its associated short boundary path.}
  \label{fig: A short path and it's associated short boundary path.}
\end{figure}

\begin{proof}
The key to this proof is corollary \ref{lemma: finding n-waveless system}.  This corollary says that it is possible to find an $n$-waveless minimal system of meridian disks $\mathbf D$ for $H$, see definition \ref{defn: n-waveless}.  There are essentially two steps to this proof. The first is to show the disk $D$ bounded by the union of a short path $\alpha$ and its associated short boundary $\alpha'\subset \bndry H$, can be homotoped into a standard position relative to an $n$-waveless minimal system of meridian disks $\mathbf D$.  The second is to show that this implies that, given a short arc $\alpha$, it is possible to construct an associated short arc $\alpha'$ from a finite list of normal curves in $\overline{\bndry H - \mathbf D}$.

As $\mathbf D$ is $n$-waveless we can assume that $\alpha$ and $\alpha'$ have minimal intersection with $\mathbf D$. Let $D'$ be a disk and $f:D'\to M$ the map such that $f(D') = D$. We can assume that $f$ is transverse to $\mathbf D$, then $\Gamma = f^{-1}(\mathbf D)$ is a set of properly embedded simple curves. The first step is to homotop $f$ to remove simple closed curves in $\Gamma$.  This is done by the usual argument of taking an innermost simple closed curve and then by irreducibility of $H$, there is a homotopy of $f$ to remove the intersection. Note that the boundary of $D$ has not moved. Thus by repeating this process all simple closed curves in $\Gamma$ can be removed.  To save on notation, the pull back of $\alpha$ and $\alpha'$ to $\bndry D'$ will also be called $\alpha$ and $\alpha'$ respectively.  For $\gamma \in \Gamma$, there are three choices either $\gamma$ is parallel to $\alpha$ or $\alpha'$ or runs between them. If there are components $\Gamma$ parallel to $\alpha$ there is an innermost one $\gamma$. That is, if $\beta \subset \alpha$ is the subarc such that $\bndry \beta = \alpha \cap \gamma$ and $D''\subset D'$ is the disk bounded by $\gamma \cup \beta$ then $D''\cap \Gamma = \gamma$.  However $D''$ gives a homotopy of $f$ to reduce the number of intersection of $\alpha$ with $\mathbf D$, giving a contradiction. Thus there are no components of $\Gamma$ parallel to $\alpha$.  Similarly if $\Gamma$ has components parallel to $\alpha'$ then there must be an innermost one $\gamma$.  Now let $\beta \subset \alpha'$ be the subarc such that $\bndry \beta = \alpha' \cap \gamma$.  Then $\beta$ is a wave in $\overline{\bndry H - \mathbf D}$.  If $\beta$ is essential then as $\mathbf D$ is $n$-waveless, $\beta$ must intersect $\tc$ at least $n/2$ times, contradicting $\alpha$ being short. Thus $\beta$ must be boundary parallel, however this contradicts $\alpha'$ having minimal intersection with $\mathbf D$.  Therefore all components of $\Gamma$ must be parallel curves running between $\alpha$ and $\alpha'$, as shown in figure \ref{fig: Pull back of a short disk.}.  Therefore we can assume that $D$ is in this standard position relative to $\mathbf D$.

\begin{figure}[h]
  \begin{center}
      \includegraphics[width=3cm]{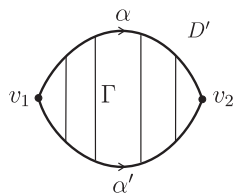}
   \end{center}
  \caption{Pull back of a short disk.}
  \label{fig: Pull back of a short disk.}
\end{figure}

Let $v_1$ and $v_2$ be the two vertices $\alpha \cap \alpha'\subset \bndry D$. From above we know that if we travel along $\alpha$ from $v_1$ to $v_2$, the order in which the disks from $\mathbf D$ are intersected is the same as if we traveled along $\alpha'$ from $v_1$ to $v_2$. As $\mathbf D$ is a waveless system for $H$, $\tc \cap \overline{\bndry H -\mathbf D}$ is a set of properly embedded curves that cut $\overline{\bndry H -\mathbf D}$ up into polygonal faces.  Therefore there are a finite number of normal curves that join any two boundary elements of $\overline{\bndry H -\mathbf D}$ and that intersect $\tc\cap\overline{\bndry H -\mathbf D}$ less than $n/2$ times.

The method for determining if a properly embedded path $\alpha$ is short, is to try and construct a short path $\alpha' \subset \bndry H$, by  building a sequence  of normal curves in $\overline{\bndry H -\mathbf D}$, that intersect $\tc$ less than $n/2$ times and that meet the disk in $\mathbf D$ in the same order as $\alpha$. Note also that when a disk of $\mathbf D$ is crossed by $\alpha$, the corresponding intersection with $\alpha'$ has the some orientation, as shown in figure \ref{fig: A short path and it's associated short boundary path.}. As there are only a finite number of such normal curves, there are only a finite number of possible sequences.
\end{proof}

\begin{figure}[h]
  \begin{center}
      \includegraphics[width=5cm]{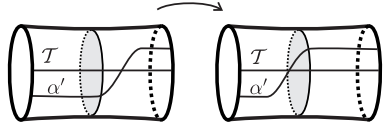}
   \end{center}
  \caption{Isotopy of an intersection across a meridian disk.}
  \label{fig: Isotopy of an intersection across a meridian disk.}
\end{figure}

It is worth noting that the normal segments used to build a short boundary path in the above proof are unique up to isotoping intersections with $\tc$ across disks in $\mathbf D$, as shown in figure \ref{fig: Isotopy of an intersection across a meridian disk.}.  This implies that the algorithm  should be reasonably efficient in detecting short paths.

Finally, we want to get a bound in terms of the boundary length on the number of faces of a disk bounded by a closed null-homotopic curve in a manifold that meets the disk-condition. Once again, let $D$ be a disk and $M$ a 3-manifold that meets the $(n_1, n_2, n_3)$ disk-condition.  Let $\mathbf n = max\{n_1, n_2, n_3\}$ and $X$ be a 2-complex that cuts $M$ up into injective handlebodies.  Let $f:D\to M$ be  a singular map and $\Gamma = \bndry D \cup f^{-1}(X)$. By lemma \ref{lemma: removing pullback graphs} we can assume that all  faces of $D$ are $(m,n)$-gons and that  $m\geq n$ for all inner faces.  Let $|\Gamma|$ be the number of faces of $\Gamma$ and $l$ be the number of vertices in $\bndry D$.

\begin{lemma}
If $\Gamma$ is as described above, then $|\Gamma| \leq 8(\mathbf n l)^2\pi $.
\end{lemma}

\begin{proof}
This lemma is proven by putting a CAT$(0)$ metric on $D$ and then using a result of Bridson and Heafliger that gives a bound on the area of CAT$(0)$ disks. The first step is to make $\Gamma$ into a graph $\Gamma'$ such that if $F$ is a face of $\Gamma'$, that is a $(m,n)$-gon, then $m \geq n$.  We already know this is true for all internal  faces.  If $F$ is a boundary $(m,n)$-gon  of $\Gamma$ such that $m < n$, we can increase its order by adding order two vertices to $F\cap \bndry D$.  Note that if the boundary faces contain $(2,n)$-gons that this process is not canonical.  However this is not a concern for the important fact is $|\Gamma|=|\Gamma'|$.  If $l$ is the number of vertices in $\bndry D \cap \Gamma$, then the number of vertices in $\bndry D \cap \Gamma'$ is at most $\mathbf n l$.

\begin{figure}[h]
  \begin{center}
      \includegraphics[width=4cm]{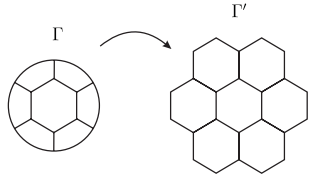}
   \end{center}
  \caption{$\Gamma$ to $\Gamma'$.}
  \label{fig:Gamma_to_Gammadash.}
\end{figure}

The next step is to put a metric on $D$ that is induced by  $\Gamma'$, however this time we need to be a little more careful. Let $F$ be a $(m,n)$-gon of $\Gamma'$.  Put a vertex $v$ at the center of $F$ and join each  vertex in $\bndry F$ to $v$ by a single edge.  Thus $F$ is cut up into $m$ triangles giving a triangulation $T$ of $D$.  Let each of these triangles be geodesic Euclidean triangles with the edges of $\bndry F$ unit length and the angle between the edges at $v$ be $2\pi/n$.  Thus all of the triangles in $F$ are isometric and the curvature is zero everywhere except at $v$, where $K(F) = 2\pi(1 - m/n)\leq 0$.  Thus $D$ is a CAT(0) metric space, as defined by Bridson and Heafliger in \cite{b&h}, and each face of $\Gamma'$ has an area of at least 1.

Bridson and Haefliger showed in ~\cite{b&h} (pg 416) that as $D$ is CAT(0) it has a triangulation $T'$ where the diameter of each face of $T'$ is at most $1$ and $|T'| \leq 8(\mathbf n l)^2$. Therefore the area of $D$ is at most $8(\mathbf n l)^2\pi$ and  $|\Gamma|\leq 8(\mathbf n l)^2\pi$.
\end{proof}

\bibliographystyle{plain}
\bibliography{bibfile}

\begin{thebibliography}{10}

\bibitem{b&h}
M.~Bridson and A.~Haefliger.
\newblock {\em Metric spaces of non-posotive curature}, volume 319 of {\em
  Grundlehren der mathematischen Wissenchaften}.
\newblock Springer, 1999.

\bibitem{C&R1}
J.~Coffey and H.~Rubinstein.
\newblock 3-manifolds built from injective handlebodies.
\newblock {\em in preparation}, 2005.

\bibitem{de1}
Max Dehn.
\newblock {\em Papers on group theory and topology}.
\newblock Springer-Verlag, New York, 1987.
\newblock Translated from the German and with introductions and an appendix by
  John Stillwell, With an appendix by Otto Schreier.

\bibitem{fhs1}
M.~Freedman, J.~Hass, and P.~Scott.
\newblock Closed geodesics on surfaces.
\newblock {\em Bull. London Math. Soc.}, 14:385--391, 1982.

\bibitem{hs1}
J.~Hass and P.~Scott.
\newblock Homotopy equivalence and homeomorphism of {$3$}-manifolds.
\newblock {\em Topology}, 31(3):493--517, 1992.

\bibitem{H&T1}
A.~Hatcher and W.~Thurston.
\newblock Incompressible surfaces in {$2$}-bridge knot complements.
\newblock {\em Invent. Math.}, 79(2):225--246, 1985.

\bibitem{Hem1}
J.~Hempel.
\newblock {\em {$3$}-{M}anifolds}.
\newblock Princeton University Press, Princeton, N. J., 1976.
\newblock Ann. of Math. Studies, No. 86.

\bibitem{Ja1}
W.~Jaco.
\newblock {\em Lectures on three-manifold topology}, volume~43 of {\em CBMS
  Regional Conference Series in Mathematics}.
\newblock American Mathematical Society, Providence, R.I., 1980.

\bibitem{Mat1}
S.~Matveev.
\newblock {\em Algorithmic topology and classification of 3-manifolds},
  volume~9 of {\em Algorithms and Computation in Mathematics}.
\newblock Springer-Verlag, Berlin, 2003.

\bibitem{rou&sa1}
C.~P. Rourke and B.~J. Sanderson.
\newblock {\em Introduction to piecewise-linear topology}.
\newblock Springer-Verlag, New York, 1972.
\newblock Ergebnisse der Mathematik und ihrer Grenzgebiete, Band 69.

\bibitem{sk1}
A.~Skinner.
\newblock The word problem in a class of non-{H}aken {$3$}-manifolds.
\newblock {\em Topology}, 33(2):215--239, 1994.

\bibitem{wa2}
F.~Waldhausen.
\newblock The word problem in fundamental groups of sufficiently large
  irreducible {$3$}-manifolds.
\newblock {\em Ann. of Math. (2)}, 88:272--280, 1968.

\end{thebibliography}

\end{document}